\documentclass[12pt]{amsart}
\usepackage{amsmath}\usepackage{amssymb}
\usepackage{amscd}
\usepackage[all]{xy}
\usepackage{color}
\usepackage[notcite,notref]{showkeys}

\newcommand{\edit}[1]{}
\def\lra{\longrightarrow}
\def\map#1{\,{\buildrel #1 \over \lra}\,}
\def\smap#1{{\buildrel #1 \over \to}}
\def\oo{\otimes}
\def\cy{{\mathrm{cy}}}
\def\ad{{\mathrm{ad}}}
\newcommand{\prim}{{\mathrm{prim}}}
\newcommand{\trf}{{\mathrm{trf}}}
\def\zar{{\mathrm{zar}}}
\def\cdh{{\mathrm{cdh}}}

\def\cE{\mathcal E}
\def\cF{\mathcal F}
\def\cO{\mathcal O}
\newcommand{\bbH}{\mathbb H}
\newcommand{\bGamma}{\bar\Gamma}
\newcommand{\bbHcdh}{\bbH_{\cdh}} 
\newcommand{\Schk}{\mathrm{Sch}/k}

\newcommand{\F}{\mathbb{F}}
\newcommand{\R}{\mathbb{R}}
\newcommand{\Q}{\mathbb{Q}}
\newcommand{\Z}{\mathbb{Z}}
\newcommand{\N}{\mathbb{N}}
\newcommand{\No}{\mathbb{N}_+}
\newcommand{\Sph}{\mathbb{S}}
   
\newcommand{\THH}{\operatorname{THH}}   
\newcommand{\TC}{\operatorname{TC}}

\def\Spec{\operatorname{Spec}}

\newcommand{\mapsp}{\mathrm{map}_{CycSp}}
\newcommand{\gh}{\mathrm{gh}}
\renewcommand{\inf}{\mathrm{inf}}

\numberwithin{equation}{section}

\theoremstyle{plain}
\newtheorem{thm}[equation]{Theorem}
\newtheorem{prop}[equation]{Proposition}
\newtheorem{lem}[equation]{Lemma}
\newtheorem{cor}[equation]{Corollary}

\newtheorem{substuff}{Remark}[equation]
\newtheorem{subcor}[equation]{Corollary}

\theoremstyle{definition}
\newtheorem{defn}[equation]{Definition}
\newtheorem{ex}[equation]{Example}

\newtheorem{Wk-mod}[equation]{$W_\Gamma(k)$-module structures}

\theoremstyle{remark}
\newtheorem{rem}[equation]{Remark}
\newtheorem{subrem}[substuff]{Remark}
\begin {document}

\title[Module structures]
{Module structure of the $K$-theory of polynomial-like rings}
\date{\today}
\author{Christian Haesemeyer}
\address{School of Mathematics and Statistics, University of Melbourne,
VIC 3010, Australia}
\thanks{Haesemeyer was supported by ARC DP-170102328}

\author{Charles A.\,Weibel}
\address{Math.\ Dept., Rutgers University, New Brunswick, NJ 08901, USA}
\email{weibel@math.rutgers.edu}\urladdr{http://math.rutgers.edu/~weibel}
\thanks{Weibel was supported by NSF grants, and the Simonyi Endowment at IAS}


\begin{abstract}
 Suppose $\Gamma$ is a submonoid of a lattice, not containing a line. In this
 note, we use the natural $\Gamma$-grading on the monoid algebra 
$R[\Gamma]$ to prove structural results about the relative $K$-theory
 $K(R[\Gamma], R)$. 
When $R$ contains a field, we prove a decomposition
 indexed by the rays
 in $\Gamma$, and a compatible action by the
 Witt vectors of $R$ for each $\N$-grading of $\Gamma$.  In
 characteristic zero, 
 there is additionally an action by Witt vectors for
%each ray of $\Gamma$. 
the  truncation set $\Gamma$.
Finally, we apply this to get a ray-like description of
$K_*(R[x_1,...,x_n])$ proposed by J.\,Davis.
\end{abstract}
\maketitle

%\edit{action of W on trucation set?\\changed abstract}
If $\Gamma$ is a commutative monoid, and $R$ is a commutative ring,
the monoid algebra $R[\Gamma]$ is naturally $\Gamma$-graded. 
In this note, we make use of this simple observation to introduce 
a natural decomposition of the relative $K$-theory groups 
$K_*(R[\Gamma],R)$ when $\Gamma$ is a normal submonoid of $\Z^n$ 
having no invertible elements except~0
and $R$ contains a field.
The decomposition is indexed
by the set of rays of $\Gamma$, where a {\it ray} is the
intersection of $\Gamma$ with a line in $\R^n=\Z^n\oo\R$.
For each ray, the summand is a continuous
module for the ring $W(R)$ of big Witt vectors of $R$.
Of course, all this generalizes the well known
$W(R)$-module structure on $K_*(R[x],R)$ 
when $\Gamma=\N=\{0,1,\dots\}$ \cite{W87}; in this case
there is only one ray.
These structures are transported from cyclic (or topological cyclic) homology
via the Chern character (cyclotomic trace) maps. 

We are motivated by the classical case where $\Gamma = \N^n$;
the cartesian decomposition in \cite{chww-bass} of the relative 
$K$-theory of $R[\N^n]=R[x_1,..,x_n]$ is compatible with the
ray-like decomposition that
%%\edit{revisit intro}
Jim Davis obtained in \cite{Davis}, using the
Farrell-Jones conjecture; this is discussed in section
\ref{sec:poly}, { which is independent of the rest of the paper}. 

Our results also apply to quotients $\bGamma=\Gamma/I$ of
$\Gamma$ by an ideal $I$. For example, if $\Gamma=\N^n$ and $I$
is generated by $\{(\cdots0,a_i,0\cdots):i=1,\dots n\}$ then we
obtain a $W(R)$-module decomposition of the $K$-theory of
$$R[\Gamma/I]=R[x_1,\dots,x_n]/(x_1^{a_1},\dots x_n^{a_n})$$
which is related to the decomposition 
of the muli-relative groups 
$$
K_*(R[\Gamma/I],(x_1,),\dots (x_n)) 
\cong W_{\Gamma/I}/J\otimes_\Z E(x_1,...,x_{n-1})
$$
found in \cite[Thm.\,1.3]{AGHL}.
{ See Propositions \ref{K-W module} and \ref{cor:decomp} below.}

We first discuss the special case where $\Q\subseteq R$ in sections
\ref{sec:HH} and \ref{sec:cdh}; the results are slightly stronger in
this case, because the ghost map is an isomorphism.

We treat the case where $\F_p\subseteq R$ in 
sections \ref{sec:TC} through 
%\ref{sec:TC-K}, \ref{Schlichtkrull} and 
\ref{sec:main}. In section \ref{sec:TC}, we 
relate the relative $K$-theory to $TC$ using the cyclotomic trace. 
In section \ref{sec:TC-K}, we
employ a result of Hesselholt and Nikolaus \cite{HN} to produce the
ray decomposition on the relative $K$-theory. Sections
\ref{Schlichtkrull} and
%\edit{quick overview of sections 3-6 -chh}
\ref{sec:main} use Schlichtkrull's description of the transfer in
topological Hochschild homology \cite{Schlicht1} to prove that the ray
decomposition is preserved by Witt vector actions arising from
$\N$-gradings of the monoid $\Gamma$.

In section \ref{sec:poly}, we extend
Davis' formula for $K_*(R[x_1,\dots,x_n])$.
In the Appendix, %\ref{sec:Gamma-Witt}, 
%Appendices \ref{sec:Gamma-Witt} and \ref{sec:Frob} develop
we develop some background material on the generalized Witt vectors 
$W_\Gamma(R)$.
% and their functorial properties.

%%%%%%%%%%%%%%%%%%%%%%%%%%%%%%%%%%%%%%%%%%%%%%%%%

\medskip
\subsection*{Acknowledgements}
We are grateful to Lars Hesselholt for his guidance on $TC$,
especially in Section \ref{sec:TC-K}.

\section{Hochschild and cyclic homology}\label{sec:HH}

In this section, we consider the Hochschild and cyclic homology 
of $A$, relative to a commutative base ring $k_0$, where
$A$ is a (possibly non-commutative) $\Gamma$-graded $k_0$-algebra.
%with $A_0=k$  commutative. 
Here $\Gamma$ is a normal submonoid of $\Z^n$
having no invertible elements except~0. 
The prototype $\Gamma$-graded algebra is the monoid ring 
$R[\Gamma]$ with $A_0=R$; 
if $I$ is an ideal of the monoid $\Gamma$, then
%\edit{edited -chh}
$R[I]$ and $R[\Gamma/I]$ are also $\Gamma$-graded algebras.
%in which $\Gamma$ is regarded as a multiplicative monoid.  

The Hochschild complex $C_*(A)$ and the kernel 
$C_*(A,A_0)$ of $C_*(A)\to C_*(A_0)$ are $\Gamma$-graded
chain complexes of $A_0$-modules.
Since $\Gamma$ has no invertible elements, the $\Gamma$-degree~0 
summand of $C_*(A)$ is the chain complex $C_*(A_0)$.
Hence their homology $HH_n(A)$ and $HH_n(A,A_0)$ 
are $\Gamma$-graded $R$-modules, $R$ being the center of $A_0$, and
$C_*(A,A_0)$ is a mixed complex of $\Gamma$-graded $k_0$-modules
with no $\Gamma$-degree~0 term.

Let $R$ be a commutative $k_0$-algebra.
Recall from Corollary \ref{idempotent}(b) that
%for any commutative ring $k$, 
$W_\Gamma(R)$ is the product
over all rays $\rho$ in $\Gamma$ of the classical Witt ring $W(R)$.
The ghost map $\gh$ sends $W_\Gamma(R)$ to the ring
$\prod_\Gamma R$ of pointed functions $\Gamma\to R$; 
see Definition \ref{ghost}. 
Thus any $\Gamma$-graded $R$-module $M$  with $M_0=0$
is a continuous $W_\Gamma(R)$-module, by
restriction of scalars along the ghost map; see Example \ref{gh-module}.

For example, when $A$ is commutative, and $A_0=R$, the 
K\"ahler differentials $\Omega^n_A=\Omega^n_{A/k_0}$
are $\Gamma$-graded with $\Omega^n_{R}$ in degree~0,
so the $\Omega^n_A/\Omega^n_{R}$ are continuous $W_\Gamma(R)$-modules.
They also form a trivial mixed complex of $\Gamma$-graded $k_0$-modules; 
see \cite[9.8.9]{WH}. 

\begin{lem}\label{ctn-module}
Let $R$ be the center of $A_0$.
For each $n$, restriction of scalars along the ghost map %$\gh$
makes each Hochschild homology group $HH_n(A,A_0)$ 
into a continuous $W_\Gamma(R)$-module,
and each cyclic homology group $HC_n(A,A_0)$
into a continuous $W_\Gamma(k_0)$-module. 
\end{lem}

\begin{proof}
This follows simply from Example \ref{gh-module}, combined with the
observation that the differential $b$ in the mixed complex
$C_*(A,A_0)$ is $R$-linear, 
\edit{was $R$-linear}
and the cyclic differential $B$ is $k_0$-linear.
\end{proof}

\begin{subrem}
If $A_0$ doesn't contain $\Q$, the ghost map is not an isomorphism,
and there are other module structures on $HH_n(A,A_0)$;
see \cite{DW} for details. Since our focus is on $K$-theory,
we do not discuss them here.
\end{subrem}

%For example, when $A$ is commutative, each K\"ahler differential module 
%$\Omega^n_A=\Omega^n_{A/R}$ is a $\Gamma$-graded $k$-module, so the
%$\Omega^n_A/\Omega^n_k$ are continuous $W_\Gamma(k)$-modules.

By Example \ref{gh-module}, %Following \cite[5.2]{DW}, 
%By Definition \ref{ghost},
for every homogeneous element $\alpha$ in $HH_n(A,A_0)$,
$HC_n(A,A_0)$ (or $\Omega^n_A/\Omega^n_{A_0}$ if $A$ is commutative)
of $\Gamma$-degree $\eta\ne0$, $\gamma\in\Gamma$ has content $c(\gamma)$ 
and $r\in k_0$, we have 
\begin{equation}\label{HHn}
r[\gamma]*\alpha = \begin{cases}
c(\gamma)\ r^e\alpha, &\mathrm{if }~ \eta=e\gamma \textrm{ in }\Gamma,
\\ 0, & \mathrm{else.}\end{cases}                                              
\end{equation}
The de\,Rham differentials $d\!:\!A/R\to\Omega^1_A/\Omega^1_R$ and 
$\Omega^n_A/\Omega^n_R\to\Omega^{n+1}_A/\Omega^{n+1}_R$ 
are homomorphisms of $W_\Gamma(k_0)$-modules,
but they are not homomorphisms of $W_\Gamma(R)$ modules
when $\Omega^1_{R/k_0}\ne0$. For example, if $t\in R$, and
$\eta=e\,\gamma\ne0$ with $\gamma$ primitive, then for $a\in A_\eta$:
\[ d(t[\gamma]*a)= d(t^e a)=t^e da+a\,d(t^e).\]
\goodbreak
%\end{subex}

\paragraph{\bf $HH$ and $HC$ of schemes.}
We can promote the above to schemes, using the fact that $HH_n(X)$
is the hypercohomology $\bbH_\zar^{-n}(X;HH)$ of $X$ with coefficients
in $HH$, the sheafification of the presheaf $C_*$ on affines; see \cite{W97}.

We now restrict to the case where $A$ is commutative, and flat over
$k_0$, with $A_0=R$ (such as { $A=R[\Gamma]$ or $R[\Gamma/I]$}),
%\edit{$A$ flat over $R$ - this is true anyway for our monoid algebras -chh}
$X$ is defined over $R$, and $R$ contains a field, writing $X\oo A$ 
for the scheme $X\times_R\,\Spec(A)$.
%$X[\Gamma]$ for the scheme $X\times_R\,\Spec(R[\Gamma])$.
By the K\"unneth formula \cite[Prop.\,9.4.1]{WH},
the shuffle product 
\[ 
HH(X)\oo HH(A,R)\map{\simeq} HH(X\oo A,X)
%HH(X)\oo_R HH(R[\Gamma],R)\map{\sim} HH(X[\Gamma],X)
\]
is a quasi-isomorphism. Taking homology yields a decomposition:
\edit{$\oo_R$ now $\oo$}
\[
HH_n(X\oo A,X)\cong
\bigoplus\nolimits_i HH_i(X)\oo HH_{n-i}(A,R).
\]
(There are no Tor-terms because $A$ is flat over $k_0$.)
\edit{flat/$k_0$}
%$X$ is defined over a field.
%cf.\ \cite[2.2]{chww-bass}.
Writing $\widetilde{HH}\oo A$ and $HH\Gamma$ for the $\Gamma$-graded sheaves
of mixed complexes associated to $HH(X\oo A,X)$ and
$HH(X\oo R[\Gamma],X)$ on $X$,
this proves: 

\begin{lem}\label{sheaf=HHA}
If $X$ is a scheme over $R$, each $HH_n(X\oo A,X)$ 
has a continuous $W_\Gamma(R)$-module structure,
natural in $X$. The same is true for the $\cdh$-hypercohomology,
and each
\[
HH_n(X\oo A,X) \to\bbH_\cdh^{-n}(X;\widetilde{HH}\oo A)
\] 
is a continuous $W_\Gamma(R)$-module map.

In particular, $\bbH_\cdh^{-n}(X;HH\Gamma)$ is a 
continuous $W_\Gamma(R)$-module. If $I$ is an ideal of $\Gamma$,
$\bbH_\cdh^{-n}(X;\widetilde{HH}\oo R[\Gamma/I])$
and $\bbH_\cdh^{-n}(X;\widetilde{HH}\oo R[I])$ are also
continuous $W_\Gamma(R)$-modules.
\end{lem}

%\begin{subrem}
%Although $HH_*(R[\Gamma/I])$ is graded by $\Gamma/I$,
%we will consider $C_*(R[\Gamma/I])$ as $\Gamma$-graded.
%\end{subrem}

Recall that if 
$K\to C\map{f}D\to TK$ is a distinguished triangle 
of chain complexes, the {\it fiber} of $f$ is $K$ and the 
mapping cone is $TK$.

\begin{cor}\label{F-version}
Let $\cF_{HH}(X)$ denote the fiber %shifted mapping cone 
of the canonical map $HH(X)\to \bbH_\cdh(X;HH)$. Then 
\edit{just $\oo$}
\[
\cF_{HH}(X\oo A/R) \cong \cF_{HH}(X) \oo HH(A,R).
%\cF_{HH}(X[\Gamma],X) \cong \cF_{HH}(X) \oo_R HH(R[\Gamma],R).
%\to\bbH_\cdh(X[\Gamma],X;HH).
\] 
Thus, its homotopy groups are also continuous $W_\Gamma(R)$-modules.
\end{cor}

\begin{rem}
If $A$ or $X$ is not flat over $k_0$, we may consider the
derived version of the Hochschild complexes $HH(A)$ and $HH(X)$.
Lemma \ref{sheaf=HHA} %\ref{ctn-module} and \ref{scheme-HH} 
and Corollary \ref{F-version}
go through in the derived setting. We leave the details to the reader.
\end{rem}

Similarly, the cyclic homology $HC_n(X)$ is the hypercohomology 
module $\bbH_\zar^{-n}(X,HC)$ of $X$ with coefficients
in the sheafification $HC$ of the complex giving cyclic homology;
see \cite{W96}. 
We write $X[\Gamma]$ for $X\oo_R R[\Gamma]$ and
%$HC\oo A/R$ (resp., $HC\Gamma$) for the kernel of 
%$HC(A)\to HC(R)$ (resp., the kernel of $HC(R[\Gamma])\to HC(R)$), and 
$\cF_{HC}(X)$ for the fiber %shifted mapping cone 
of  the canonical map $HC(X)\!\to\!\bbH_\cdh(X;HC)$.

\begin{lem}\label{scheme-HC}
If $\Q\subseteq R$ and
$X$ is a scheme over $R$, each $HC_n(X[\Gamma],X)$
has a continuous $W_\Gamma(R)$-module structure, natural in $X$. 
The same is true for the $\cdh$-hypercohomology, and each
\[
HC_n(X[\Gamma],X) \to\bbH_\cdh^{-n}(X;HC[\Gamma])
%HC_n(X[\bGamma],X) \to\bbH_\cdh^{-n}(X,HC\bGamma)
\]
is a continuous $W_{\Gamma}(R)$-module map. In addition,
%\edit{check!\\checked 09/04}
the homotopy groups of $\cF_{HC}(X[\Gamma],X)$
% $\cF_{HC}[\bGamma],X)$ 
are also continuous $W_\Gamma(R)$-modules.

If $I$ is an ideal of $\Gamma$ then the homotopy groups of 
$\cF_{HC}(X[I],X)$ and $\cF_{HC}(X[\Gamma/I],X)$ are also
continuous $W_\Gamma(R)$-modules.
\end{lem}

\begin{proof}
Because $\Q[\Gamma]$ has at least one $\N$-grading with $\Q$ in degree~0,
the $S$ operator acts as zero on $HC_*(A,\Q)$; see \cite[9.9.1]{WH}.
This means that $HC_*(A,\Q)$ is a trivial $HC_*(\Q)$-comodule.
By \cite[(3.2)]{Kassel}, this implies that
 $HC_*(R[\Gamma],R)\cong HH_*(R)\oo HC_*(\Q[\Gamma],\Q)$,
%\edit{explained ok? 8/28\\ yes, couple typos corrected 09/04}
and similarly with $R$ replaced by $X$. 
Since a continuous $W(R)$-module is a colimit of $W(R)/V_mW(R)$-modules,
the tensor product of an $R$-module with a continuous 
$W(\Q)$-module is a continuous $W(R)$-module.
Hence  $HC_*(R[\Gamma],R)$
has the structure of a continuous $W_\Gamma(R)$-module.
\end{proof}

\begin{rem}
Suppose that $k$ is a field containing $\Q$, and that 
$X$ is a toric variety over $k$. 
%\edit{use ``a field $k$''\\and not $R$?\\ yes, changed 09/04}
By \cite[2.4]{chww-vorst}, the Zariski sheaf $HH$
%Let $HH$ denote the sheaf of Hochschild chain complexes.
%by definition \cite{W97}, $HH_n(X) = \bbH_\zar^{-n}(X,HH)$.
is the product of complexes $HH^{(i)}$\!, $i\ge0$, and 
$\bbH_\cdh(X;HH^{(i)}) \simeq \bbH_\cdh(X;\Omega^i)[i]$.
As pointed out in \cite[2.1--2]{chww-vorst},
%$(C_*)_\cdh\simeq \prod \Omega_\cdh^i[i]$ in the derived category,
%and hence 
$\bbH_\cdh^n(X;HH)$ is the direct sum of the groups
$H_\cdh^{n+i}(X;\Omega^i)$. 
Similarly, $HC$ is the product of complexes $HC^{(i)}$\!, $i\ge0$, and
$\bbH_\cdh(X;HC)$ is the product of the 
complexes $\bbH_\cdh(X;HC^{(i)})$.
\end{rem}

%Let $HH\Gamma$ denote the presheaf of complexes $X\mapsto HH(X[\Gamma],X)$.
%The shuffle product $HH(X)\oo HH(k[\Gamma],k)\to HH(X[\Gamma],X)$,
%which is a quasi-isomorphism by the proof
%of Lemma \ref{scheme-HH}, sheafifies with respect to $X$ to 
%yield a shuffle product $HH_\cdh(X)\oo HH(k[\Gamma],k)\to HH\Gamma_\cdh(X)$.
%Taking hypercohomology, we have a quasi-isomorphism
%\[
%  \bbH_\cdh(X,HH) \oo HH(k[\Gamma],k)\map{\sim} \bbH_\cdh(X,HH\Gamma).
%\]
%It follows that the groups
%$\bbH_\cdh^*(X,HH\Gamma)$ are $W_\Gamma(k)$-modules, and that
%$HH_n(X[\Gamma],X) \to\bbH_\cdh^*(X,HH\Gamma)$ are $W_\Gamma(k)$-module maps.

%\smallskip
%We now sheafify with respect to the $\cdh$-topology. By definition,
%$\bbH^*_\cdh(X,HH)$ is the hypercohomology of the sheafification of 
%the Hochschild complex $C_*$ for the $\cdh$-topology, 
%or equivalently, of a complex  $\bbH_\cdh(X,HH)$ 
%of $\cdh$-sheaves quasi-isomorphic to $(C_*)_\cdh$.
%Similarly, $\bbH^*_\cdh(X,HC)$ is the hypercohomology of a complex 
%$\bbH_\cdh(X,HC)$.

\goodbreak
%\newpage
\section{$K$-groups in characteristic~0}\label{sec:cdh}

Let $\Schk$ be the category of schemes of finite type over a field $k$ 
%\edit{using $k$ over $k_0$}
of characteristic~0, and write $X[\bGamma]$ for 
$X\times\Spec(\Z[\bGamma])$, where 
$\bGamma$ is $\Gamma$, $I$ or $\Gamma/I$.
The purpose of this section is to show that the relative $K$-groups 
\[  
K_n(X[\bGamma],X)=K_n(X[\bGamma])/K_n(X)
%K_n(X\times T,X)=K_n(X\times_k T_k)/K_n(X)  
\]
are continuous $W_{\Gamma}(k)$-modules, by relating them to 
cyclic homology relative to $k_0=\Q$.

Let $\cF\!_K(X)$ denote the fiber of the canonical map 
$K(X)\!\to\!\bbH_\cdh{\!}(X,{\!}K)$. 
By \cite[Thm.\,1.6]{chww-vorst}, the Chern character $K\to HN$ induces
a weak equivalences $\cF_K(X)\smap{\sim}\vert\cF_{HC}(X)[1]\vert$
for each $X$ in $\Schk$, where $|C|$ denotes the associated spectrum
of a chain complex $C$. 
Given Lemma \ref{scheme-HC}, %was Corollary \ref{cor:better},
the following result is a rephrasing of \cite[Thm.\,5.6]{chww-toric}; 
in that theorem, $\bGamma$ is assumed to be finitely generated as a monoid
but, since $K$-theory commutes with direct limits, this is not a problem.

\begin{prop}\label{K-W module}
For all $n$, 
$K_n(X[\bGamma],X) \cong \bbH_\cdh^{1-n}(X,\cF_{HC\bGamma})$.
Hence $K_n(X[\bGamma],X)$ is a continuous $W_{\bGamma}(k)$-module.

In particular, $K_n(X[\bGamma],X)$ decomposes as a direct sum
of $W(k)$-modules $K_n(X[\bGamma],X)_\rho$, 
indexed by the rays $\rho$ of $\Gamma$.
\end{prop}

\begin{proof}
By \cite{Haes}, $\bbHcdh(X;K)\simeq KH(X)$, which is homotopy invariant.
Because $k[\bGamma]$ has an $\N$-grading with $k$ in degree~0,
$\bbHcdh(X;K)\simeq\bbHcdh(X[\bGamma];K)$, i.e.,
$\bbHcdh(X[\bGamma],X;K)\simeq~0.$
Hence 
\[
|\cF_{HC}(X[\bGamma],X)[1]|
\stackrel{\simeq}{\leftarrow}
\cF_K(X[\bGamma],X) \map{\simeq} K(X[\bGamma],X). %\qedhere
\]
The final assertion follows from Lemma \ref{oplus}.
\end{proof}

\begin{ex}
If $\bGamma=\N^n$, this proves that $K(k[x_1,...,x_n])/K_n(k)$
is a continuous $W_{\Gamma}(k)$-module.
We give more detail in Corollary \ref{cor:NnK} below.
\end{ex}

%For the computation in the next example, compare
%\cite[Proposition 5.6]{chww-toric}
%\edit{our toric paper}
%\edit{There is no Prop 5.6 in \cite{chww-monoids}\\ also checked
% \cite{chww-toric}, but does not fit -chh\\maybe [6,5.6]? }
%module structure on the relative $K$-theory of a toric variety
% $X[\bGamma]$ in characteristic~0.

\begin{ex}\label{ex:xy=zz}
Let $\Gamma$ be the submonoid of $\N^2$ generated by 
$x=(1,0)$, $y=(1,2)$ and $z=(1,1)$ so that $k[\Gamma]$ is
$k[x,y,z]/(z^2=xy)$.  By \cite[4.3]{chww-cones},
%\edit{our cones paper}
$K_1(k[\Gamma],k)\cong k$
with the $k$ in $\Gamma$-degree $(2,2)$;
$a\in k$ corresponds to the class of the matrix
$\bigl( \begin{smallmatrix} 1+az &ax\\ay &1-az
        \end{smallmatrix} \bigr)$.
%$\left(\binom{1+az}{ay}\binom{ax}{1-az} \right)$.
More generally, if $i\,>\,0$ then $K_{2i}(\Q[\Gamma],\Q)=0$  and
$K_{2i-1}(\Q[\Gamma],\Q)\cong\Q$.

Since $k$ is a field, it follows that %$K_1(k[\Gamma],k)\cong k$ and,
$K_q(k[\Gamma],k)\cong \oplus_{r\ge0}\,\Omega_k^{q-1-2r}.$
Replacing $k$ by $k[\N^p]=k[x_1,...,x_p]$, it follows that 
\[
K_1(k[\Gamma\wedge\N^p],k[\N^p])\cong k[\N^p]
\cong K_1(k[\Gamma],k) \oplus \bigoplus\nolimits_{\rho\subset\N^p}\ tk[t]
\] 
(the rays $\rho$ in $\N^p$ correspond to maximal inclusions
of $k[t]$ into $k[\N^p]$), 
%$(1,1)\times\N^p$),
and more generally that
$K_q(k[\Gamma\wedge\N^p],k)$ is a sum 
$\oplus_{r\ge0}\,\Omega_{k[x_1,...,x_p]}^{q-1-2r}$, where
\[
\Omega_{k[x_1,...,x_p]}^{q-1-2r} \cong \bigoplus_{i+j=q-1-2r}
\Omega^i_k\otimes \Omega^j_{\Q[x_1,...,t_p]}.
\]
\end{ex}

%\newpage
\section{Topological cyclic homology}\label{sec:TC}
%{Topological cyclic homology, and 
%relative $K$-theory in arbitrary characteristic.}

The purpose of this section is to relate the $K$-groups
$K_n(X[\bGamma],X)$ to topological cyclic homology,
in preparation for Sections \ref{sec:TC-K}
and \ref{sec:main}.

We begin by recalling the definition of $\TC(R)$, following
Nikolaus--Scholze \cite[II.1]{NS}. 
Let $T$ denote the circle group.
A {\it cyclotomic spectrum}
is a spectrum $E$ together with $T$-equivariant maps
$E\to E^{tC_p}$ for all primes $p$. %($T$ is the circle.) 

For example, if $R$ is a ring then
$\THH(R)$ is a cyclotomic spectrum, and  
$\TC(R)=\TC(\THH(R))$ is defined in \cite[II.1.8]{NS}
to be the mapping spectrum
\[ TC(R) = \mapsp(\Sph,\THH(R)).\]

We may also define $\TC^-(R) =\THH(R)^{hT}$; defining $TP(R)$ as
the product over all primes of $\left(\THH(R)^{tC_p}\right)^{hT}$,
there is a functorial fiber sequence (see \cite[II.1.9]{NS}):
\begin{equation}\label{SBI}
\TC(R) \to \TC^-(R) \to TP(R).
\end{equation}

Let $X$ be a finite dimensional noetherian scheme.
There is no difficulty in defining $\THH(X)$ by Zariski descent,
and we define $\TC(X)$ %=\TC(\THH(X))$ 
to be $\mapsp(\Sph,\THH(X))$.

The cyclotomic trace map $K\to \TC$ of \cite[1.2]{NS} induces
a homotopy cartesian square:
\begin{equation}\label{square:K-TC}
\xymatrix{
K(X)\ar[d] \ar[r]& \bbHcdh(X,K) \ar[d] \\
\TC(X) \ar[r]     & \bbHcdh(X,\TC).
}\end{equation}
%is a homotopy cartesian square. 
This follows from Land--Tamme \cite[B.3]{LT},
given the fact \cite[3.6]{LT} that 
the fiber $K^{\inf}$ of the cyclotomic map
$K\to TC$ is ``truncating'';  a
spectrum-valued functor $E$ on perfect $\infty$-categories
is called {\it truncating} if for each connective $A_\infty$-spectrum $A$,
$E(A)\to E(\pi_0A)$ is an equivalence,
and $E$ sends exact sequences to fiber sequences.
Thus the homotopy fibers $\cF_K$ and $\cF_{TC}$ of the
horizontal maps in \eqref{square:K-TC} are equivalent.
(This was proven under restrictive hypotheses
by Geisser and Hesselholt in \cite[Thm.\,C]{GH}.)

By \cite[6.3]{KST} or \cite[A.5]{LT}, $\bbHcdh(X,K)\cong KH(X)$
for any finite dimensional noetherian scheme $X$. 
By homotopy invariance, if $\bGamma$ is either 
$\Gamma$, $I$ or $\Gamma/I$ then
$KH(X[\bGamma],X)=0$. 
This proves:
 
\begin{prop}\label{K-TC}
For any finite-dimensional noetherian scheme $X$,
$$K(X[\bGamma],X) \map{\simeq} \cF_{TC}(X[\bGamma],X).$$
%$K(k[\bGamma],k) \map{\simeq} \cF_{TC}(k[\bGamma],k)$.
\end{prop}
\goodbreak

\section{Ray decomposition of $K$-groups in characteristic $p$.}
\label{sec:TC-K}

Let $\Z$ denote the spectrum $H\Z$, regarded as a cyclotomic spectrum
with trivial action by the circle group $T$.
%there is a map  of cyclotomic spectra, $\Z\to \THH(X)$.
As pointed out in \cite{HN}, when $k=\Z/p$,
there is a $T$-equivariant map %from $\Z$ %(with the trivial $T$-action) 
\[
\Z \to \Z/p\simeq \tau_{\ge0}TC(k)\to THH(k).
\]
(This fails for $k=\Z$.)
By the argument of \cite[Sec.\,1]{HN}, this gives an 
equivariant decomposition of cyclotomic spectra 
for each scheme $X$ over $\Z/p$ :
\[
\THH(X[\bGamma]) \cong \THH(X)\oo_\Sph B^{\cy}\bGamma_+
\cong \THH(X) \oo_\Z (\Z\oo_\Sph B^{\cy}\bGamma_+).
\]
Here $B^{\cy}\bGamma$ is the cyclic nerve of $\bGamma$; 
it is also called $N^{\cy}(\bGamma)$.
% and $\oo_\Sph$ is the symmetric monoidal product of spectra
With respect to this decomposition,
each Frobenius $\varphi:\THH\to \THH^{tC_p}$ factors 
as a composition $\varphi\oo\tilde\varphi$,
where $\varphi(X):\THH(X)\to \THH(X)^{tC_p}$, and
$\tilde\varphi: B^{\cy}\bGamma\to (B^{\cy}\bGamma)^{hC_p}$
is the unstable Frobenius.

By \cite[Sec.\,1]{HN}, there is an equivalence between 
the $\infty$-category of chain complexes with $T$-action,
$D(\Z)^{BT}$, and the $\infty$-category of mixed complexes over $\Z$.
%\edit{removed sphere after re-reading \cite{HN},
% this $\otimes$ means ``free abelian on" -chh}
Under this equivalence, 
the object $\Z\oo B^{\cy}\bGamma_+$ of $D(\Z)^{BT}$
corresponds to the mixed complex $C_*(\Z[\bGamma])$ of $\Z[\bGamma]$.
The $\bGamma$-grading gives a decomposition of
$\THH(X[\bGamma],X)$ into a sum of terms 
$\THH(X)\oo C_*(X[\bGamma],X)_\rho$ indexed by 
the rays $\rho$ of $\Gamma$.

\begin{prop}\label{TC-decomp}
For every scheme $X$ over $\Z/p$, each of the spectra
$TC(X[\bGamma],X)$, $TC^-(X[\bGamma],X)$
and $TP(X[\bGamma],X)$ decomposes as the sum of factors
$TC(X[\bGamma],X)_\rho$, etc.,
indexed by the rays $\rho$ of $\Gamma$.
There are similar decompositions for  $TC^-$ and $TP$.
%The same is true for $TC^-$ and $TP$.
\end{prop}

\begin{proof}
The complex $C_*(\Z[\bGamma])$ is $\bGamma$-graded, and each unstable
Frobenius $\tilde\varphi$ preserves the rays $\rho$ of $\bGamma$; 
for each $\eta\in\bGamma$, $\tilde\varphi$ sends the
degree $\eta$ summand to the degree $\eta^p$ summand.
Hence $TC^-(X[\bGamma]) \map{\varphi^{hT}} TP(X[\bGamma])$ 
decomposes as a sum of maps $TC^-_\rho\to TP_\rho$.
It follows that $TC(\Z[\bGamma],\Z)$ is the sum over all rays $\rho$
of $TC(\Z[\bGamma])_\rho$.
\end{proof}

The ray-indexed decompositions in Proposition \ref{TC-decomp}
are compatible with the cdh fibrant replacement
functor $TC\to\bbH_\cdh(-,TC)$. Thus they give a
decomposition of the fiber $\cF_{TC}(X[\bGamma],X)$
as the sum over all rays $\rho$ in $\Gamma$
of $\cF_{TC}(X[\bGamma],X)_\rho$. 
%Here $X$ is $\Spec(\Z[\bGamma])$.
%Consequently, the fiber $\cF_{TC}$ ...

Combining this with Proposition \ref{K-TC}, this proves

\begin{cor}\label{cor:decomp}
For any finite-dimensional noetherian scheme $X$ over $\Z/p$, 
the spectrum $K(X[\bGamma],X)$ decomposes as a sum of spectra
$K(X[\Gamma],X)_\rho$, 
indexed by the rays $\rho$ of $\Gamma$.
Hence there are natural decompositions
$K_n(X[\bGamma],X)\!\cong\!\oplus_\rho K_n(X[\bGamma],X)_\rho$.
%the groups $K_n(X[\bGamma],X)$ decompose as a direct sum
%$\oplus_\rho K_n(X[\bGamma],X)_\rho$.
%indexed by the rays $\rho$ of $\Gamma$.
\end{cor}

%\newpage %%%%%%%%%%%%%%%%%%%%%%%%%%%%%%%%%%
\section{The transfer in $THH$}\label{Schlichtkrull}

Let $G$ be a pointed, abelian monoid, and $K$ a submonoid such that,
as a $K$-set, $G$ is free of rank $m$. Schlichtkrull \cite{Schlicht1}
constructed a transfer map $THH(RG)\to THH(RH)$ in the special case when
$G$ and $K$ are groups (with a disjoint basepoint); the purpose
of this section is to show that much of his construction goes through
in the monoid setting.

\begin{defn}\label{def:EGG}
Let $\cE G$ denote the simplicial set $EG\times_G G^{\ad}$,
where $G^{\ad}$ denotes $G$ as a trivial $G$-set. It is the
sum $\bigvee_{g\in G} EG\times\{g\}$.
\end{defn}

In \cite[Sec.\,6]{Schlicht1}, Schlichtkrull constructs maps
\[ 
\cE G \map{\trf} \Gamma^+ (EG\times_K G^{\ad}), \qquad
EG\times_K G^{\ad}  \map{\kappa} \Gamma^+\cE K.
\]
The first is the topological transfer map of the $m$-sheeted covering space
$EG\times_K G^{\ad} \to \cE G$, and the second is the
$K$-set projection $G\cong \bigvee_i K\gamma_i\to K\gamma_1\cong K$,
where $\{\gamma_1,...,\gamma_m\}$ is a $K$-basis of the $K$-set $G$.
He writes $\widetilde{\mathrm{Res}}$ for the composition
\[
\Gamma^+(\kappa)\circ\trf:
\Sigma^\infty(\cE G) \to \Gamma^+\Sigma^\infty(\cE K).
\]
We will apply his results to $G=\bGamma\times\langle t\rangle$ 
and $K=\bGamma\times\langle t^m\rangle$.

The following result is based on Schlichtkrull's calculations
in \cite{Schlicht1}. 
For simplicity, we write $\bGamma[t]$ for the monoid 
$\bGamma\times\langle t\rangle$ and fix  $\gamma t^n$ in $\bGamma[t]$.
Recall from \cite[2.2]{chww-char=p} that $THH(R[\bGamma])$ 
is $S^1$-equivariantly $THH(R)\wedge B^\cy\Gamma$,
and that the $\bGamma$-grading comes from 
$B^\cy\bGamma = \bigvee_{g\in\bGamma} B^\cy(\bGamma;g)$.
There is a similar equivalence for $B^\cy(\bGamma[t])$.
%See \cite[2.2]{chww-char=p}.

\begin{prop}\label{prop:trf}
For $\trf:THH(R)\wedge B^\cy\bGamma[t]\smap{}THH(R)\wedge B^\cy\bGamma[t]$:\\
%$\gamma t^m$ in $\bGamma[t]$, acting on $THH(R)\wedge B^\cy\bGamma[t]$:\\
(a) if $m\,\not\vert\, n$ then $\trf=0$ on the $\gamma t^n$ component
of $THH(R)\wedge B^\cy\bGamma[t]$;\\
(b) if $m\vert n$, $\trf$ sends the $\gamma t^n$-component
of $THH(R)\wedge B^\cy\bGamma[t]$ to itself.
\end{prop}

\begin{proof}
If $n=0$, we get $B^\cy\bGamma[t]_\gamma\simeq 
B^\cy\bGamma_\gamma\wedge B^\cy\langle t\rangle_1$
and $B^\cy\langle t\rangle_1 \simeq S^0$ 
by \cite[2.9.2]{chww-char=p} and \cite{Schlicht1}.
The same is true for $B^\cy\langle t^m\rangle$, 
and the transfer is simply the identity.

If $n>0$, we get $B^\cy\bGamma[t]_{\gamma t^n} = 
B^\cy\bGamma_\gamma\wedge B^\cy\langle t\rangle_{t^n}$;
by \cite[3.20--21]{Rognes} this is 
$B^\cy\bGamma_\gamma\wedge B^\cy\langle t,1/t\rangle_{t^n}$.
In addition, the following diagram commutes.
\begin{equation*}
\xymatrix{
THH(R)\wedge B^\cy\bGamma_\gamma\wedge B^\cy\langle t\rangle_{t^n} 
\ar[r]^\simeq\ar[d]^\simeq & 
THH(R)\wedge B^\cy\bGamma_\gamma\wedge B^\cy\langle t,1/t\rangle_{t^n}\ar[d]^\simeq
\\
THH(R[t])_{t^n}\wedge B^\cy\bGamma_\gamma  \ar[r]^\simeq\ar[d]^{\trf\wedge 1} &
THH(R[t,1/t])_{t^n}\wedge B^\cy\bGamma_\gamma\ar[d]^{\trf\wedge 1}
\\ 
THH(R[t^m])\wedge B^\cy\bGamma \ar[r]^{\simeq}& THH(R[t^m,1/t^m])\wedge B^\cy\bGamma 
}\end{equation*}
The vertical compositions are the transfer maps
from $THH(R\bGamma[t])_{\gamma t^n}$ to $THH(R\bGamma[t^m])$ and
$THH(R\bGamma[t^m,1/t^m])$, respectively.
%$THH(R[\bGamma[t])_{\gamma t^m}\to THH(R[\bGamma[t^n])$ and
%$THH(R[\bGamma[t])_{\gamma t^m}\to THH(R[\bGamma[t^n,1/t^n])$.
Now apply \cite[Thm.\,A]{Schlicht1} to see that the lower right
vertical map is 0 unless $m\vert n$, and maps to the 
$\gamma t^n$-component otherwise.
\end{proof}

\begin{subrem}
In \cite[Prop.\,6.1]{Schlicht1}, Schlichtkrull assumes that
$G$ and $K$ are groups 
(with disjoint basepoints), so that $\cE G\simeq BG_+$, 
and compares $\widetilde{\mathrm{Res}}$
%$\bGamma^+(\kappa)\circ\trf:\cE G \to \cE K$
to the restriction map $N^\cy G\to N^\cy K$ using a map
$\phi:N^\cy G \to \cE G$.
This approach does not carry over to monoids.
\end{subrem}
%Showing that $W(R)$-action preserves the $\bGamma$-grading
%boils down to %checking  that the composite 

\begin{cor}\label{THH-W}
Let $\phi$ and $\psi$ be graded $R\bGamma$--algebra endomorphisms 
of $R\bGamma [t]$. The $\bGamma$-grading on $THH(R\bGamma)$ is 
%\edit{added $\phi$ and $\psi$ for $[r]$ and $V_m$ later -chh}
preserved by the composition
\[
THH(R\bGamma) \map{i} THH(R\bGamma[t]) \map{\psi\circ\trf\circ\phi} 
THH(R\bGamma[t])    \map{t=1} THH(R\bGamma).
\]
%preserves the $\bGamma$-grading on $THH(R\bGamma)$.
%\edit{``Do it for THH''}
\end{cor}

\begin{cor}\label{TC-W}
Let $\phi$ and $\psi$ be graded $R\bGamma$--algebra endomorphisms 
of $R\bGamma[t]$. The composition
\[
TC(R\bGamma) \map{i} TC(R\bGamma[t]) \map{\psi\circ\trf\circ\phi}
TC(R\bGamma[t])    \map{t=1} TC(R\bGamma)
\]
%preserves the $\bGamma$-grading on $THH(R\bGamma)$.
%transfer $\trf$  on $TC(R[\bGamma],R)$
preserves the ray decomposition on $TC(R\bGamma)$. 
\end{cor}

\begin{proof}
This follows because $TC$ is constructed from the 
%\edit{``Do it for $TC$''}
cyclotomic structure on $THH$ in a way that preserves the ray decomposition
(but {\it not} the $\bGamma$-grading). Indeed, the transfer %Frobenius 
maps the $\gamma$-component of $TC(R[\bGamma],R)$
to the $\gamma^p$ component; see the discussion
on page 144 of \cite{HN}.
\end{proof}

\begin{cor}\label{K-W}
Let $\phi$ and $\psi$ be graded $R\bGamma$--algebra 
endomorphisms of $R\bGamma [t]$. The ray decomposition 
of $K(R\bGamma)$ is preserved by the composition
\[
K(R\bGamma) \map{i} K(R\bGamma[t]) \map{\psi\circ\trf\circ\phi}
K(R\bGamma[t])    \map{t=1} K(R\bGamma).
\]
%transfer $\trf$  on $K(R[\bGamma],R)$
\end{cor}

\begin{proof} 
Since everything is natural in $R$--algebras, 
the transfer $\trf$ and composition $\psi\circ\trf\circ\phi$ 
act on the $\bGamma$-graded spectrum
%\edit{``do it for $K$''}
%there is a $W(R)$--action
$\bbH_\cdh(\Spec R, TC(\cO[\bGamma],\cO))$, 
and the action is compatible with the ray decomposition.
Hence the same is true for the homotopy fiber $\cF(R[\bGamma],R)$ 
of the map
\[
TC(R\bGamma,R)\to \bbH(\Spec R,TC(\cO[\bGamma],\cO)).
\]
Finally, we observe that the isomorphism
$K_*(R\bGamma,R) \cong \pi_* \cF_{TC}(R\bGamma,R)$
is compatible with the $\bGamma$-decomposition and with the
actions of $\trf$ and $\psi\circ\trf\circ\phi$.
\end{proof}

%\medskip\hrule\medskip  %%%%%%%%%%%%%%%%%%%%%%%%%%%%%%%%%%%%%
\section{Ray components are $W(k)$-modules}\label{sec:main}

In \cite{W87}, a functorial formula was given for a
continuous action of the Witt vectors $W(k)$ %$(1-rt^m)$ 
on the relative $K$-theory $K_*(A,A_0)$ of any $\N$-graded 
$k$-algebra $A$. 
This applies to $\Gamma$-graded $A$ once we choose a 
homomorphism $p:\Gamma\to\N$ such that $p(\gamma)\ne0$ 
for all $\gamma$ with $A_\gamma\neq 0$; %with $p(\Gamma)\ne0$;
these exist by our assumption that $\Gamma$ contains no line.
Each grading $\Gamma\to\N$ also defines a $W(k)$-action
on $HH_*(A,A_0)$, etc.; see \cite{DW}.
Let $i:A\to A[t]$ denote the injection sending $a\in A_n$ to $at^n$.
Recall from \cite{W87} that
multiplication by $\omega\in W(k)$ on $K(A,A_0)$ is the composition
\begin{equation}\label{eq:K(A)}
K(A,A_0) \map{i} K(A[t],A) \map{\omega} K(A[t],A) \map{t=1} K(A,A_0).
\end{equation}

Here is one way to choose a grading $\Gamma\to\N$.
Recall that $\Gamma$ is the monoid of lattice points in the 
``dual cone'' of a rational polyhedral cone $\sigma\subset\R^n$.
Set $\Gamma^\vee=\{w\in\Z^n: w\cdot\Gamma\ge0\}$, and
choose  $w\in\Gamma^{\vee}$ with 
$w\cdot v>0$ for all $0\neq v\in\Gamma$.
%Given a ray $\rho$ in $\Gamma$, with primitive vector $v$,
%choose a $w$ in $\sigma$ to minimize $w\cdot v$ subject to
%%$w\cdot v>0$; 
This gives an $\N$-grading on any $\Gamma$-graded 
algebra $A$ such that $A_0$ is the ring of 
elements of ($\N$--) degree $0$, and hence
a $W(k)$-module stucture on any $K_*(A,A_0)$, and in particular on
$K_*(k[\Gamma],k)$.
However, as pointed out in \cite[5.3]{W87}, the $W(k)$-module structure
can depend on the choice of $w$. 

Let $R$ be a commmutative $k$-algebra, 
%\edit{added 9/6}
and $X$ a finite-dimensional noetherian scheme with $R=H^0(X,\cO_X)$.
Note that $W(R)$ is a $W(k)$-algebra.
Fixing a choice of $w$, we obtain:

\begin{thm}\label{thm:main}
Each $K_n(R[\bGamma],R)$ and $K_n(X[\bGamma],X)$ is a 
continuous $W(R)$-module; the latter decomposes as a sum 
$\oplus_\rho K_n(X[\bGamma],X)_\rho$
of continuous $W(R)$-modules,
%\edit{does $W(k)$ act or $W(R)$ -chh}
indexed by the rays $\rho$ of $\Gamma$.
%\edit{tweaked 2/26,3/6}
%For each ray $\rho$ in $\Gamma$, the summand 
%$K(X[\Gamma],X)_\rho$ is a continuous $W(k)$-module.
\end{thm}

\begin{proof}
Since every $\omega\in W(R)$ is a product 
$\omega=\prod_m(1-r_mt^m)$, it suffices
to consider the action of each $\omega=(1-rt^m)$
on $K_n(R[\bGamma],R)$; the case $K_n(X[\bGamma],X)$ is similar.
Set $A=R[\bGamma]$. 
%and observe that in  \eqref{eq:K(A)}, 
%$i$ sends rays of $\Gamma$ to rays of $\Gamma[t]$.
%Thus it suffices to consider the action of $(1-rt^m)$ on $K_n(A[t],A)$.

By \cite{W87},
$\omega$ acts on $K(A[t],A)$ as $F_m[r]V_m$. %see \cite{W87}. 
Now $[r]$ and $V_m$ correspond to the graded $A$-algebra homomorphisms 
%
%\edit{$A$-algebra hom's emphasised -chh}
%
$t\mapsto rt$ and $t\mapsto t^m$, 
%they clearly preserve the ray decomposition. 
and the operator $F_m$
corresponds to the transfer map $K(A[t])\to K(A[t^m])$
followed by the equivalence associated with $A[t^m]\cong A[t]$.
%
%\edit{this seems fine now 09/04}
%\edit{Proof re-written because I couldn't see why $K(A[t],A)$ 
%decomposes wrt rays, please CHECK -chh}
%Finally, the transfer 
%preserves the ray decomposition of $K(A[t],t)$
%
Thus, the composition (\ref{eq:K(A)}) with $\omega = (1 - rt^m)$ 
preserves the ray decomposition by Corollary \ref{K-W}.
\end{proof}

\bigskip\goodbreak
%\newpage
\section{Polynomial rings}\label{sec:poly}

{\noindent The results in this section are independent of the rest of this paper.}
{Let $k$ be any unital ring.}
According to the Fundamental
Theorem of $K$-theory \cite[III.2.2, V.8.2]{WK},
$K_q(k[t_1,...,t_n])/K_q(k)$ is the direct sum (for $1\le i\le n$) of
$\binom{n}{i}$ copies of $N^iK_{q}(k)$.
%This holds for every unital ring $k$.
%Similarly, $K_q(k[t_1,t_1^{-1},...,t_n,t_n^{-1}])$ is the sum of
%$\binom{n}{i}$ copies of $K_{q-i}(k)$ ($0\le i\le n$) and 
%$\binom{n}{i}\binom{n-i}{j}2^{i}$ copies of $N^iK_{q-j}(k)$.

Recall that $k[\Z^n] = k[t_1,t_1^{-1},...,t_n,t_n^{-1}]$ is
the group ring of the free abelian group $\Z^n$,
and the Farrell--Jones conjecture in $K$-theory is known for the 
free abelian group $\Z^n$; see Quinn \cite{Quinn}. Using this,
J.\ Davis proved in \cite{Davis} (compare \cite{LR}) that 
$K_q(k[\Z^n])$ is the direct sum of 
$\binom{n}{r}$ copies of $K_{q-r}(k)$ ($r=0,...,n$)
and the direct sum over the rays $\rho\subset\Z^n$ of
many copies of $NK_{q-r}(k)$ ($r=0,...,n-1$):
\begin{equation*}
\bigoplus\nolimits_{r=0}^n \binom{n}{r}K_{q-r}(k) \quad\oplus\quad
\bigoplus\nolimits_{\rho \subset\Z^n} 
\bigoplus\nolimits_{r=0}^{n-1} 
\binom{n-1}{r} NK_{q-r}(k),
\end{equation*}
(Every maximal cyclic subgroup of $\Z^n$ is the union of two rays.)
Using the notation that { $L^rK_q=K_{q-r}$ and} 
$L^rNK_q = NK_{q-r}$, this becomes
\begin{equation}\label{eq:Davis}
K_q(k[\Z^n]) \cong (1+L)^nK_q(k) \oplus
\bigoplus\nolimits_{\rho\subset\Z^n} (1+L)^{n-1}NK_q(k).
\end{equation}
When $n=1$ there are only 2 rays, so we recover
the usual formula for $K_q(k[x,1/x])$.
The formula for $n=2$ implies the formula:
\[
K_q(k[x,y],k)\cong\bigoplus\nolimits_{\rho\subset\N^2} 
NK_q(k)\oplus NK_{q-1}(k).
\]
Let us write $\No^p$ for the { sub-semigroup} %monoid 
of elements $(n_1,...,n_p)$ in $\N^p$ with all $n_i>0$.
Since $\N^2-\{0\}=(\No\times0) \cup (0\times\N_+) \cup \No^2$, 
the above formula implies that 
\[
N^2K_q(k) \cong 
\bigoplus\nolimits_{\rho\subset\No^2} NK_q(k)\oplus NK_{q-1}(k).
\]
%
%Let $\N^r_+$ denote the set of all $(x_1,...,x_r,0,...0)$
%in $\N^n$ with $x_i>0$ for $i=0,...,r$.

\begin{thm}\label{thm:NnK}
For every unital ring $k$,
$$N^{n}K_q(k) \cong \bigoplus_{\rho\subset \No^n} (1+L)^{n-1}NK_q(k).$$
\end{thm}

This was conjectured by Davis in \cite{Davis}, 
who proved it %and proven by Davis
in \cite[Cor.\,15]{Davis} when the $NK_i$ are countable torsion groups.
We first give a quick argument to
establish Theorem \ref{thm:NnK} when $k$ is a $\Q$-algebra.

\begin{proof}[Proof of \ref{thm:NnK} for commutative $\Q$-algebras]
Identifying $t\Q[t]$ and $\Omega^1_{\Q[t]/\Q}$, the formula of
\cite[4.2]{chww-bass} is that for $p>0$:
%Identifying $tk[t]$ and $\Omega^1_{k[t]/k})$, this becomes:
\[
N^{p+1}K_q(k) \cong x_1\cdots x_p\Q[x_1,...,x_p]\oo_{\Q} (1+L)^pNK_q.
\]
Plugging in $NK_q=TK_q\oo t\Q[t]$ (from \cite[0.1 and 4.2]{chww-bass})
yields Theorem \ref{thm:NnK} in this case.
\end{proof}

%When $k$ is a $\Q$-algebra, a different-looking formula
%was given in \cite[4.2]{chww-bass}; for example, 
%$N^2K_q(k)$ is a direct sum of copies of
%$NK_q(k)$, indexed by monomials $t^n$ ($n\ge1$), and
%copies of $NK_{q-1}(k)$, indexed by $t^{n-1}dt$ ($n\ge1$),
%respectively. This agrees with Davis' formula if we use
%the formula $NK_q(k)\cong TK_q\otimes tk[t]$ of \cite[(0.3)]{chww-bass}.
%(Use $tk[t]^{\oo2}\cong xyk[x,y]$.)
%
%More generally, we have
%\[
%N^{p+1}K_q(k) \cong \bigoplus_{r=0}^p NK_{q-r}(k)\otimes\wedge^r k^p
%\otimes (tk[t])^{\otimes p-r}\otimes(\Omega^1_{k[t]/k})^{\otimes r}.
%\]
%Identifying $tk[t]$ and $\Omega^1_{k[t]/k})$, this becomes:
%\[
%N^{p+1}K_q(k) \cong x_1\cdots x_pk[x_1,...,x_p]\oo (1+L)^pNK_q.
%\]
%Plugging in $NK_q=TK_q\oo tk[t]$ yields Theorem \ref{thm:NK} in this case.

For the proof of \ref{thm:NnK} for arbitrary $k$,
recall that the wreath group $W_n=(\Z/2)\wr S_n$ of order $2^n n!$ acts 
as automorphisms on the monomials in $k[x_1,x_1^{-1},...,x_n,x_n^{-1}]$;
the subgroup $(\Z/2)^n$ interchanges the $x_i$ and $x_i^{-1}$
pairwise, and the subgroup $S_n$ permutes the $x_i$. It also acts on $\Z^n$
(by signed permutation matrices), on the primitive vectors in $\Z^n$,
 and on the rays in $\Z^n$.
For each $r<n$, let $R_r$ denote the set of all rays in
$(x_1,...,x_r,0,...,0)$ with all $x_i$ positive. The orbit of $R_r$
under $W_n$ has $2^r{\binom nr}$ disjoint components, and
the stabilizer of $R_r$ is $S_r\times W_{n-r}$.
{ For example, when $r=1$ and $n=2$, the orbit of $(1,0)\in R_1$
consists of the 4 points $(\pm1,0)$ and $(0,\pm1)$.
}

Similarly, the group $S_n$ acts on $\N^n$, and the orbit of $R_r$
has ${\binom nr}$ disjoint components, and the stabilizer of $R_r$ is
$S_r\times S_{n-r}$.

\begin{proof}[Proof of Theorem \ref{thm:NnK}]
We proceed by induction on $n$, the case $n=1$ being trivial.
Let $Q$ denote the right hand side of Theorem \ref{thm:NnK}.

Recall that the $W_{r}$-orbit of $\No^r$ in $\Z^r$ has order $2^r$.
By induction, for $r<n$ the images of the $2^r$ embeddings of
$N^rK_q(k) \subset K_q(k[\N^r])$ into $K_q(k[\Z^r])$ give embeddings of 
{ $(1+L)^{r-1}NK_q(k)$}  %$(1+L)^{n-r}N^rK_q(k)$
into $K_q(k[\Z^n])$. The image
is the sum of copies of $(1+L)^{n-1}NK_q(k)$, indexed by
the rays in the $2^r$ copies of $R_r$.
%primitive vectors in $\N^r_+$, 
As $r$ varies, we get a subgroup of 
$K_q(k[\Z^n])$ whose quotient is $2^r$ copies of $Q$.
%the right hand side of Theorem \ref{thm:NnK}.

This induces a $W_n$-equivariant isomorphism
from $\oplus_{2^n} N^nK_q(k)$ to $\oplus_{2^n} Q$.
%$N^{n}K_q(k) \cong \oplus_{\rho\in \N^n} (1+L)^{n-1}NK_q(k)$.
Taking invariants yields a map from $N^nK_q(k)$ to
\[
{(\oplus_{2^n}Q)^{W_n}=}
\oplus_{\rho\in \No^n} (1+L)^{n-1}NK_q(k).
\qedhere
%K_q(k)\oplus \bigoplus\nolimits_{\N^r_+} NK_q
\]
\end{proof}

{
\begin{rem}
The referee points out that the case of general $n$ 
follows from the case $n=2$ by induction.
When $n=2$, the rays on $(1,0)$ and $(0,1)$ correspond to the
two canonical inclusions of $NK_q(k)$ in $K_q(k[x,y])$,
given by the subrings $k[x]$ and $k[y]$ of $k[x,y]$.
\end{rem}
}

\smallskip
\begin{cor}\label{cor:NnK}
For every unital ring $k$, 
%$K_q(k[x_1,...,x_n])$ is isomorphic to 
\[
K_q(k[x_1,...,x_n]) \cong
K_q(k)\oplus \bigoplus\nolimits_{r=1}^n
\bigoplus\nolimits_{\rho\subseteq \No^r} {\binom nr}(1+L)^{r-1}NK_q(k).
\]
\end{cor}

\noindent
For example, the summmand for $r=1$ is $n$ copies of $NK_q(k)$,
as $\rho=\No^1$.
When $n=2$ we have the decomposition: 
\begin{align*}
 K_q(k[x,y])\cong & K_q(k)\oplus 2NK_q(k)\oplus 
\bigoplus\nolimits_{\rho\subset\No^2}(NK_q(k)\oplus NK_{q-1}(k))
\\
\cong & K_q(k)\oplus \bigoplus\nolimits_{\rho\subset\N^2} NK_q(k) ~\oplus~ 
\bigoplus\nolimits_{\rho\subset\No^2}NK_{q-1}(k).
\end{align*}
This is a strengthening of Bass' formula
\[
K_q(k[x,y])\cong (1+N)^2K_q(k) = K_q(k)\oplus 2NK_q(k)\oplus N^2K_{q}(k).   
\]

\bigskip %\newpage
\appendix%%%%%%%%%%%%%%%%%%%%%%%%%%%%%%%%%%%%
\section{Witt vectors on affine monoids}\label{sec:Gamma-Witt}

Let $\Gamma$ be a submonoid of $\Z^n$ having no invertible elements 
except 0. %and set $\Gamma=\Gamma\cup\{0\}$.
We assume that $\Gamma$ is {\it normal} in the sense that
if $v\in\Z^n$ and $dv\in \Gamma$ for some positive integer $d$ then
$v\in \Gamma$. The following definition is generalized from \cite{Borger}, where $\Gamma = \N^r$.
{ 
\begin{defn}  A subset $S$ of $\Gamma\setminus\{0\}$ is a
$\Gamma$-{\it truncation set} if, for all natural numbers $d$ and $s\in \Gamma$,
$ds\in S$ implies $s\in S$. 
\end{defn}

\begin{subrem}
	Let $I$ be an ideal in $\Gamma$. Then the subset $\Gamma \setminus (I \cup \{0\})$ is a $\Gamma$-truncation set. Indeed, if $\gamma\in\Gamma$ and $d\gamma\notin I$ for some $d\in\N$, then $\gamma\notin I$. 
\end{subrem}

}
In this section, we fix a commutative ring $k$ and
define the rings $W_\Gamma(k)$ { and $W_{\Gamma/I}(k)$}, 
so that when $\Gamma=\N$ { and $I = \{n\in \N \; \vert \; n\geq d\}$} we recover the classical ring $W(k)$ of big Witt vectors 
{ and truncated Witt vectors} on $k$.
Much of this section is
adopted from \cite{AGHL}, \cite{HH} and \cite{Wick}. 
% as $\Gamma-\{0\}$ is an example of a {\it truncation set}.

%\medskip
\begin{defn}\label{def:prim}
A nonzero element $v$ of $\Gamma$ is called {\it primitive} if it is not
a multiple of any other element of $\Gamma$.
Each primitive element $v$
generates a maximal cyclic submonoid $\rho$ of $\Gamma$, called a {\em ray}.
Because $\Gamma$ has no invertible elements,
each ray has the form $\N v$ for a unique primitive $v$. %$v\in\Gamma$;
As a pointed set, $\Gamma$ is the wedge of its rays.

Fix a commutative ring $k$. The ring of {\it big Witt vectors} 
$W_\Gamma(k)$ over $k$ %form a ring, whose
has for its underlying set $k^\Gamma$, the pointed functions 
$\Gamma\map{a} k$ (functions which send $0$ to $0$). 
Addition and multiplication in $W_\Gamma(k)$ are determined 
by Lemma \ref{W-ring} below, exactly as in the classical ring $W(k)$.
We write $a_\gamma$ for $a(\gamma)$, so  $a$ can be
written as the vector $\{a_\gamma\}_{\gamma\in\Gamma}$. 
{ If $S\subset\Gamma \setminus \{0\}$ is a
  $\Gamma$-truncation set, we can similarly define the truncated Witt
  vectors $W_S(k)$; the underlying set is the set of functions $S\to
  k$. If $S = \Gamma \setminus (I \cup \{0\})$ for an ideal $I$, then
  we write $W_{\Gamma/I}(k)$ for the ring $W_S(k)$. In this case,
  $S\cong \Gamma/I \setminus \{0\}$.}
\end{defn}
%\paragraph{\bf Rays.}
%The maximal cyclic submonoids of $\Gamma$ have the form 
%$\R\gamma\cap \Gamma$ (for some $\gamma\in\Gamma$) and are called
%the {\it rays} of $\Gamma$; as a pointed set $\Gamma$ is 
%the wedge of its rays. 

Here are some special elements of $W_\Gamma(k)$. For each $\gamma\in\Gamma$,
we write $[\gamma]$ for the Kronecker $\delta$-function sending 
$\gamma$ to~1 and $\eta$ to 0 if $\eta\ne\gamma$, 
For $r\in k$, $r[\gamma]$ sends $\gamma$ to $r$ and $\eta\ne\gamma$ to~0.
We write $\delta_{\prim}$ for the Kronecker $\delta$-function sending 
primitive elements to~1 and non-primitives to 0.

For any nonzero $\gamma\in \Gamma$, we define its {\it content} $c(\gamma)$ 
to be the least common multiple of $\{ e\in\N: \gamma\in e\cdot\Gamma\}$.
As an element of $\Z^n$, $\gamma=(e_1,...,e_n)$, and
$c(\gamma)=\gcd\{e_i\}$.
It is easy to see that $\gamma$ is primitive if and only if $c(\gamma)=1$.
Because $\Gamma$ is normal, $v=\gamma/c(\gamma)$ is a primitive element
of $\Gamma$ and the ray $\rho=\R\gamma\cap \Gamma$ is $\N v$.

%formal power series $\sum_\gamma a_\gamma[\gamma]$.
%$\sum_\gamma a_\gamma\delta_\gamma$.

{ Given a $\Gamma$-truncation set $S$, let $\prod_S k$ denote the set $k^{S}$ of functions, 
made into a ring 
under the termwise sum $(a+b)_\gamma = a_\gamma+ b_\gamma$ and product 
$(ab)_\gamma=a_\gamma b_\gamma$. This ring is just the product of 
copies of the ring $k$, indexed by $S$.}

\begin{defn}\label{ghost}
The {\it ghost map} { $\gh:W_S(k)\to \prod_S k$} is the 
function taking $a=\{a_\gamma\}$ to the function 
{ $\gh(a):S\to k$} defined by
\[
\gh(a): \eta \mapsto \sum_{e\gamma=\eta} c(\gamma) a_\gamma^e.
\]
The $\eta$-component of $\gh(r[\gamma])$ is 
$c(\gamma)r^e$ if $\eta=e\gamma$, and 0 otherwise.
\end{defn}
%for the sum in $k^\Gamma$ of $[v]$ for primitive $v$; 
\goodbreak

%\begin{subrem}
%The classical ghost map $(1+tk[[t]])^\times\to \prod_1^\infty k$
%is the homomorphism sending $f(t)$ to the coefficient
%sequence of $tf'(t)/f(t)$; $f(t)=(1-r_et^e)$ is sent to the function 
%$\N\to k$ sending $n$ to $e\,r_e^d$ if $n=de$, and~0 otherwise.
%
%When $\Gamma=\N$, there is a bijection 
%$(1+tk[[t]])^\times\to \prod_1^\infty k,$ 
%sending the formal product 
%$\prod_e(1-r_et^e)$ to the sequence $(r_1,r_2,...)$.
%Composing the ghost map \ref{ghost} with this bijection yields
%the classical ghost map, which is the homomorphism 
%$(1+tk[[t]])^\times\to \prod_1^\infty k$ sending $f(t)$ to the coefficient
%sequence of $tf'(t)/f(t)$; $f(t)=(1-r_et^e)$ is sent to the function 
%$\N\to k$ sending $n$ to $e\,r_e^d$ if $n=de$, and~0 otherwise.
%\end{subrem}
%on the group  $(1+tk[[t]])^\times$, via the bijection %between
%$(1+tk[[t]])^\times\to \prod_1^\infty k$ sending the formal product 
%$\prod_e(1-r_et^e)$ to the sequence $(r_1,r_2,...)$.
%Composing this bijection with the ghost map \ref{ghost} yields the 
%classical ghost map $f\mapsto tf'(t)/f(t)$; this homomorphism sends 
%$(1-r_et^e)$ to the function $\N\to k$ sending $n$ to $e\,r_e^d$
%if $n=de$, and~0 otherwise.

\begin{lem}\label{W-ring}
\cite[2.3]{AGHL}
There is a unique way, functorial in $k$, 
to put a ring structure on the set { $W_S(k)$}
so that { $\gh:W_S(k)\to \prod_S k$} is a ring homomorphism.

The identity element of the ring { $W_S(k)$} is $\delta_\prim$,
and $\gh(\delta_\prim)$ is the identity of { $\prod_S k$.}
\end{lem}

\begin{proof}
{ $W_S(k)$ and $k^S$ are the products over all rays 
$\rho=\N v$ in $\Gamma$
of the sets $W_{\N v\cap S}(k)$ and $k^{\N v\cap S}$, respectively.
Since the ghost map respects this decomposition, we may assume that
$S$ is a truncation set in $\N$. In this case, the result is classical.}
\end{proof}

\begin{subcor}\label{idempotent}
a) 
If $v$ is primitive, $[v]$ is an idempotent element of $W_\Gamma(k)$;
the idempotent $\gh([v])$ of $\prod_{\Gamma\setminus\{0\}} k$ is the
projection onto $\prod_{\N v\setminus \{0\}} k$, and $[v]\cdot W_\Gamma(k)\cong W(k)$.

\noindent b) 
$W_\Gamma(k)$ is the product over all rays $\rho=\N v$ in $\Gamma$
of the classical ring $W(k)$ of big Witt vectors:
\[ 
W_\Gamma(k) = \prod_{\N v}W_{\N v}(k)\cong 
\prod_{\mathrm{rays} \rho} W(k). 
\]
In particular, when  $\Gamma=\N$, the ring $W_{\N}(k)$ is
the classical ring $W(k)$ of big Witt vectors over $k$. 
\end{subcor}

%\begin{ex}
When $\Q\subseteq k$, the ghost map is an isomorphism:
$W_\Gamma(k)\cong \prod_{\Gamma\setminus\{0\}} k$.
This is because it is an isomorphism in the classical setting.
%\end{ex}

\begin{subrem}
By abuse of notation, for each primitive $v\in\Gamma$,
we will write $t_v$ for $[v]$ and $t_v^n$ for $[nv]$, so that the
abelian group structure on $W_{\N v}(k)$ may be interpreted as
the multiplicative group $(1+t_vk[[t_v]])^\times$.
Since the abelian group structure on $W_\Gamma(k)$ corresponds to
the group $\prod_{v}(1+t_vk[[t_v]])^\times$, 
and not a subgroup of the units in $k[[\Gamma]]$,
the power series interpretation of $W_\Gamma(k)$ 
is inconvenient in our setting.
\end{subrem}

\begin{subrem}
If $I\subset\Gamma$ is an ideal of $\Gamma$, let $\Gamma/I$ denote 
the quotient monoid { with underlying set} $\Gamma\setminus I$.
The quotient map $\Gamma\to\Gamma/I$ yields a 
ring homomorphism $W_\Gamma(k)\to W_{\Gamma/I}(k)$ sending 
the vector $\{r_\gamma\}_{\gamma\in\Gamma\setminus\{0\}}$
to $\{r_\gamma\}_{\gamma\not\in I\cup\{0\}}$. 
%In the terminology of Borger
%\cite{Borger}, \cite{HH}, \cite{Ang} and \cite{AGHL}, the semigroup
%$\Gamma/I-\{0\}$ is a {\it truncation set}: if $\gamma\in\Gamma$ and 
%$d\gamma\not\in I$ then $\gamma\not\in I$.
\end{subrem}

\paragraph{\bf Continuous $W_\Gamma(k)$-modules}
We say that a $W_\Gamma(k)$-module $A$ is {\it continuous} if
for each element $a\in A$ there is a finite subset $S_a$ of $\Gamma$
such that the products $[\gamma]*a$ vanish 
for all $\gamma\not\in S_a$.

Since $W_\Gamma(k)$ is a product ring, every $W_\Gamma(k)$-module $A$
is contained in a product $\prod A_\rho$ of $W(k)$-modules, indexed by the
rays $\rho=\N v$ in $\Gamma$;
%primitive elements $v\in\Gamma$; 
$A_\rho=[v]\cdot A$.
The following characterization is immediate.

\begin{lem}\label{oplus}
A $W_\Gamma(k)$-module $A$ is a continuous module if and only if
$A=\oplus_\rho A_\rho$ and each $A_\rho$ is a continuous $W(k)$-module.
\end{lem}

\begin{ex}\label{gh-module}
%Set $\Gamma'=\Gamma-\{0\}$.
Any $\Gamma$-graded $k$-module 
$A=\oplus_{\gamma\ne0}A_\gamma$ may be regarded
as a continuous  module over the product ring $\prod_{\Gamma} k$, 
with $(x_\gamma)\in \prod_\Gamma k$ acting via multiplication 
by $x_\gamma$ on the component $A_\gamma$, $\gamma\ne0$.
The restriction along the ghost map $W_\Gamma(k)\to \prod_\Gamma k$
defines a %(continuous)
$W_\Gamma(k)$-module structure on $A$. 
By Definition \ref{ghost}, if $a\in A_\eta$ then
\[
r[\gamma]*a = \begin{cases}
c(\gamma) r^e a, &\mathrm{if }~ \eta=e\gamma \textrm{ in }\Gamma,
\\ 0, & \mathrm{else.}\end{cases}                                              
\]
Since $\Gamma$ is a submonoid of $\mathbb{Z}^n$,
this module structure is continuous. 

If $\Q\subseteq k$, the category of continuous 
$W_\Gamma(k)$-modules is equivalent to the category of 
$\Gamma$-graded $k$-modules $A=\oplus_{\gamma\ne0} A_\gamma$, where
$\{ A_\gamma\}_{\gamma\in\Gamma}$ is a family of $k$-modules.
%This follows from the fact
%that $\gh\!:W_\Gamma(k)\cong\prod_{\gamma\in\Gamma}k$, 
%every continuous $W_\Gamma(k)$-module has the form
\end{ex}

%_____________________________________


\begin{thebibliography}{10}

\bibitem{Ang}
V. Angeltveit, %,Vigleik(5-ANU-MI)
Witt vectors and truncation posets,
{\em Theory Appl. Categ.} 32 (2017), paper 8, 258–285. 

\bibitem{AGHL}
V. Angeltveit, T. Gerhardt, M. Hill and A. Lindenstrauss,
On the algebraic $K$-theory of truncated polynomial algebras
in several variables,
{\em J. $K$-theory} 13 (2014), 57--81. arXiv:1206.0247

%\bibitem{Bloch}
%S. Bloch,
%Algebraic $K$-theory and Crystalline Cohomology,
%{\em Publ. Math. I.H.E.S.} 47 (1978), 188--268.


\bibitem{Borger}
J. Borger, 
The basic geometry of Witt vectors, I: The affine case,
{\em Algebra Number Theory} 5 (2011), 231--285. 

%\bibitem{C67}
%P. Cartier,
%Groupes formels associ\'es aux anneaux de Witt g\'en\'eralis\'es,
%C. R. Acad. Sci. Paris 265 (1967) S\'er.A, 49--52.

\bibitem{chww-vorst}
G. Corti\~nas, C. Haesemeyer, M. Walker and C. Weibel,
$K$-regularity, $cdh$-fibrant Hochschild homology, and a conjecture of Vorst, 
{\em J. AMS} {\bf21} (2008), 547--561.

\bibitem{chww-bass}
G. Corti\~nas, C. Haesemeyer, M. Walker and C. Weibel,
Bass' $NK$ groups and $cdh$-fibrant Hochschild homology,
{\em Invent. Math.} 181 (2010), 421--448.

\bibitem{chww-toric}
G. Corti\~nas, C. Haesemeyer, M. Walker and C. Weibel,
The $K$\!--theory of toric varieties,
{\em Trans.\ AMS} 361 (2009), 3325--3341.

\bibitem{chww-char=p}
G. Corti\~nas, C. Haesemeyer, M. Walker and C. Weibel,
The $K$\!--theory of toric varieties in positive characteristic $p$,
{\it J. Topology} 7 (2014), 247-286.

\bibitem{chww-cones}
G. Corti\~nas, C. Haesemeyer, M. Walker and C. Weibel,
The $K$\!--theory of cones of smooth varieties, 
{\em J. Alg. Geom.} 22 (2013), 13--34. 

\bibitem{chww-monoids}
G.\,Corti\~nas,  C.\ Haesemeyer, M.\,Walker and C.\,Weibel,
Toric Varieties, Monoid Schemes and $cdh$ descent,
{\it J. reine angew. Math.} 698 (2015), 1--54.


%\bibitem{Danilov}
%V. Danilov,
%\newblock The geometry of toric varieties,
%{\em Russian Math. Surveys} 33 (1978), 97--154.

\bibitem{Davis}
J.F. Davis,
Some remarks on Nil groups in algebraic K-theory,
available at arXiv:0803.1641 (2008).


\bibitem{DW}
B. Dayton and C. Weibel,
Module structures on the Hochschild and Cyclic Homology of Graded Rings,
pp. 63--90 in {\it Algebraic $K$\!--theory and Algebraic Topology}, 
NATO ASI Series C, {\bf407}, Kluwer, 1993.

\bibitem{Fu}
W. Fulton,
{\it Introduction to Toric Varieties}, 
Princeton Univ. Press, 1993.

\bibitem{GH}
T. Geisser and L. Hessselholt,
On the vanishing of negative $K$-groups,
Math. Ann. 348 (2010), 707--736.


%\bibitem{Gu05}
%J. Gubeladze,
%K-theory of toric varieties revisited,
%{\em J. Homotopy Rel. Structures} 9 (2014), 9--23.

%\newblock The nilpotence conjecture in $K$-theory of toric varieties,
%{\em Inventiones Math.} 160 (2005), 173--216.

\bibitem{Haes}
C. Haesemeyer,
Descent properties of homotopy $K$-theory,
{\em Duke Math. J.} 125 (2004), 589--619.

\bibitem{HH}
L. Hesselholt, 
The big de Rham--Witt complex,
{\em Acta Math.} 214 (2015), 135--207.

%\bibitem{HH-cusp}
%L. Hesselholt, 
%planar cuspidal curves

\bibitem{HN}
L. Hesselholt and Nikolaus,
Algebraic $K$-theory of planar cuspidal curves,
 K-theory in algebra, analysis and topology,
AMS Contemp. Math. 749 (2020),  139–148. 

\bibitem{Kassel}
C. Kassel,
Cyclic homology, comodules, and mixed complexes,
{\em J. Alg.} 107 (1987), 195--216.

\bibitem{KST}
M. Kerz, F. Strunk and G. Tamme,
Algebraic K-theory and descent for blow-ups,
Invent. Math. 211 (2018), 523--577.

\bibitem{LR}
W. L\"uck and D. Rosenthal,
On the $K$- and $L$-theory of hyperbolic and 
virtually finitely generated groups,
{\em Firum Math.} 26 (2014), 1565--1609. 

\bibitem{LT}
M. Land and G. Tamme,
On the $K$-theory of pullbacks,
Annals of Math. 190 (2019), 877--930.

\bibitem{NS}
T. Nikolaus and P. Scholze,
On topological cyclic homology,
Acta Math. 221 (218), 203---409.

\bibitem{Quinn}
F. Quinn, 
Algebraic K-theory over virtually abelian groups,
{\em J. Pure Appl. Algebra} 216 (2012), 170--183.
%was: Hyperelementary assembly for K-theory of virtually abelian groups

\bibitem{Rognes}
J. Rognes,
\newblock Topological logarithmic structures,
\newblock In {\em New topological contexts for Galois theory and 
algebraic geometry (BIRS 2008)}, volume~16 of
{\em Geom. Topol. Monogr.}, pages 401--544, Geom. Topol. Publ,
Coventry,  2009.


\bibitem{Schlicht1}
C. Schlichtkrull,
The transfer map in topological Hochschild homology,
{\em J. Pure Appl. Algebra} 133 (1998), 289--316.


%\bibitem{Schlicht2}
%C. Schlichtkrull,
%Transfer maps and the cyclotomic trace,
%{\em Math. Annalen} 336 (2006), 191--238.

 
\bibitem{WH}
C. Weibel,
{\it An introduction to homological algebra},
Cambridge U. Press, 1994.

%\bibitem{Wnil} 
%C. Weibel,
%Nil K-theory maps to Cyclic Homology,
%{\it Trans. AMS} 303 (1987), 541-558.

\bibitem{W87}
C. Weibel,
Module structures on the $K$-theory of graded rings,
{\em J. Algebra} 105 (1987), 465--483.

\bibitem{W96}
C. Weibel,
Cyclic homology for schemes,
{\it Proc. AMS} 124 (1996), 1655--1662.

\bibitem{W97}
C. Weibel,
The Hodge filtration and cyclic homology,
{\it $K$\!--theory} 12 (1997), 145--164.

\bibitem{WK}
C. Weibel,
{\it The $K$-book}, %an introduction to algebraic $K$-theory},
Grad. Studies in Math. 145, AMS, 2013.

\bibitem{Wick}
K. Wickelgren,
Cartier's first theorem for Witt vectors on $\Z^n_{\geq 0} - 0$,
pp. 321--328 in {\it Algebraic topology: applications and new directions},
AMS Contemp. Math. 620, 2014.

\end{thebibliography}
\end{document}